\documentclass[graybox]{svmult}
\usepackage{mathptmx}       
\usepackage{helvet}         
\usepackage{courier}        
\usepackage{makeidx}         
\usepackage{graphicx}        
\usepackage{multicol}        
\usepackage[bottom]{footmisc}
\usepackage{amsmath,amssymb,amsfonts}        

\usepackage{latexsym,a4,bm,bbm}

\makeindex             

\def\bbbc{{\mathbb C}}
\def\bbbq{{\mathbb Q}}

\def\bbbz{{\mathbb Z}}
\def\bbbf{{\mathbb F}}
\def\fqstar{\bbbf_q^{\times}}
\def\q{{\mathbbm{q}}}

\def\boldalpha{\bm{\alpha}}
\def\boldbeta{\bm{\beta}}

\def\mod#1{({\rm mod\ }#1)}

\def\dim{{\rm dim}}

\def\is{\equiv}
\def\mod#1{({\rm mod}\ #1)}
\def\floor#1{\lfloor#1\rfloor}
\def\tr{{\rm Tr}}
\def\galQ{{\rm Gal}(\overline{\bbbq}/\bbbq)}
\def\galK{{\rm Gal}(\overline{\bbbq}/K)}
\newenvironment{proof1}[1]{\noindent {\bf Proof #1: }}{\hfill\break\medskip\hfill$\Box$\medskip}
\newenvironment{example1}{\smallskip \refstepcounter{subsection}
{\bf Example \thesubsection.}}
{\hfill\break\smallskip\hfill$\diamondsuit$\smallskip}
\newtheorem{theorem1}[subsection]{Theorem}
\newtheorem{definition1}[subsection]{Definition}
\newtheorem{lemma1}[subsection]{Lemma}
\newtheorem{proposition1}[subsection]{Proposition}
\newtheorem{assumption}[subsection]{Assumption}

\begin{document}
\title*{Fields of definition of finite hypergeometric functions}
\author{Frits Beukers}
\institute{University of Utrecht, \email{f.beukers@uu.nl}}
%
%
\maketitle

\abstract{Finite hypergeometric functions are functions of a finite field
$\bbbf_q$ to $\bbbc$. They arise as Fourier expansions of
certain twisted exponential sums and were introduced independently
by John Greene and Nick Katz in the 1980's. They have many properties
in common with their analytic counterparts, the hypergeometric functions.
One restriction in the definition of finite hypergeometric functions
is that the hypergeometric parameters must be rational numbers whose
denominators divide $q-1$. In this note we use the symmetry in the
hypergeometric parameters and an extension of the exponential sums to
circumvent this problem as much as posssible.
}

\section{Introduction}
In the 1980's John Greene \cite{greene} and Nick Katz \cite{katz}
independently introduced functions from finite fields to the complex
numbers which can be interpreted as finite sum analogues of the
classical one variable hypergeometric functions. These functions,
also known as Clausen-Thomae functions are determined by two
multisets of $d$ entries in $\bbbq$ each. We denote them by 
$\boldalpha=(\alpha_1,\ldots,\alpha_d)$ and $\boldbeta=(\beta_1,\ldots,\beta_d)$.
Throughout we assume that these sets have empty intersection when considered
modulo $\bbbz$. The Clausen-Thomae functions satisfy a linear differential
equation of order $d$ with rational function coefficients. See \cite{beukersheckman}.

Let $\bbbf_q$ be the finite field with $q$ elements.
Let $\zeta_p$ be a primitive $p$-th root of unity and define the additive
character $\psi_q(x)=\zeta_p^{\tr(x)}$ where $\tr$ is trace from $\bbbf_q$ to $\bbbf_p$.
For any multiplicative character
$\chi:\fqstar\to\bbbc^{\times}$ we define the Gauss sum
$$g(\chi)=\sum_{x\in\fqstar}\chi(x)\psi_q(x)\;.$$
Let $\omega$ be a generator of the character group on $\fqstar$.
We use the notation
$g(m)=g(\omega^m)$ for any $m\in\bbbz$. Note that $g(m)$ is periodic in $m$ with
period $q-1$. Note that the dependence of $g(m)$ on $\zeta_p$ and $\omega$ is
not made explicit. 
Very often we shall need characters on $\fqstar$ of a given order.
For that we use the notation $\q=q-1$ so that a character of order $d$ can be given by
$\omega^{\q/d}$ for example, provided that $d$ divides $\q$ of course.

Now we define finite hypergeometric sums. Let again $\boldalpha$ and $\boldbeta$ be multisets
of $d$ rational numbers each, and disjoint modulo $\bbbz$. We need the following crucial
assumption.

\begin{assumption}\label{divisionassumption}
Suppose that
$$(q-1)\alpha_i,(q-1)\beta_j\in\bbbz$$
for all $i$ and $j$.
\end{assumption}

\begin{definition1}[Finite hypergeometric sum]\label{definitionHq}
Keep the above notation and Assumption \ref{divisionassumption}.
We define for any $t\in\bbbf_q$,
$$H_q(\boldalpha,\boldbeta|t)={1\over 1-q}\sum_{m=0}^{q-2}
\prod_{i=1}^d\left({g(m+\alpha_i\q)g(-m-\beta_i\q)
\over g(\alpha_i\q)g(-\beta_i\q)}\right)\ \omega((-1)^dt)^m\;.$$
\end{definition1}

It is an exercise to show that the values of $H_q(\boldalpha,\boldbeta|t)$ are
independent of the choice of $\zeta_p$. 

The hypergeometric sums above were considered without the normalizing factor
$(\prod_{i=1}^dg(\alpha_i\q)g(-\beta_i\q))^{-1}$
by Katz in \cite[p258]{katz}. Greene, in \cite{greene},
has a definition involving Jacobi sums which, after some
elaboration, amounts to
$$\omega(-1)^{|\boldbeta|\q}q^{-d}\prod_{i=1}^d{g(\alpha_i\q)g(-\beta_i\q)\over g(\alpha_i\q-\beta_i\q)}
\ H_q(\boldalpha,\boldbeta|t)\;,$$
where $|\boldbeta|=\beta_1+\cdots+\beta_d$.
The normalization we adopt in this paper coincides with that of Dermot McCarthy,
\cite[Def 3.2]{mccarthy}.

Let
$$A(x)=\prod_{j=1}^d(x-e^{2\pi i\alpha_j}),\qquad B(x)=\prod_{j=1}^d(x-e^{2\pi i\beta_j}).$$
An important special case is when $A(x),B(x)\in\bbbz[x]$. In that case we say that the
hypergeometric sum is defined over $\bbbq$. Another way of describing this case is
that $k\boldalpha\is\boldalpha\mod{\bbbz}$ and $k\boldbeta\is\boldbeta\mod{\bbbz}$
for all integers $k$ relatively prime to the common denominator of the $\alpha_i,\beta_j$.
In other words, multiplication by $k$ of the $\alpha_i\mod{\bbbz}$ simply permutes these
elements. Similarly for the $\beta_j$. From work of Levelt \cite[Thm 3.5]{beukersheckman}
it follows that in such a case the monodromy group of the classical hypergeometric equation
can be defined over $\bbbz$. It also turns out that hypergeometric sums defined over $\bbbq$
occur in point counts in $\bbbf_q$ of certain algebraic varieties, see \cite[Thm 1.5]{BCM}
and the references therein. It is an easy exercise to show that
$H_q(\boldalpha,\boldbeta|t)$ is independent of the choice of $\omega$ (it is already
independent of the choice of $\psi_q$).

One of the obstacles in the definition of finite hypergeometric sums over $\bbbq$
is Assumption \ref{divisionassumption} which has to be made on $q$,
whereas one has the impression that such sums can be defined for any $q$ relatively
prime with the common denominator of the $\alpha_i,\beta_j$. This is resolved in
\cite[Thm 1.3]{BCM} by an extension of the definition of hypergeometric sum.
The idea is to apply the theorem of Hasse-Davenport to the products of Gauss sums
which occur in the coefficients of the hypergeometric sum.
Another way of dealing with this problem is given by McCarthy, who uses the 
Gross-Koblitz theorem which expresses Gauss sums as values of the $p$-adic
$\Gamma$-function.

\begin{theorem1}[Gross-Koblitz]
Let $\omega$ be the inverse of the Teichm\"uller character. Let $\pi^{p-1}=-p$
and $\zeta_p$ such that $\zeta_p\is1+\pi\mod{\pi^2}$.
Let $\Gamma_p$ be the $p$-adic Morita
$\Gamma$-function. Let $q=p^f$ and $g_q(m)$ denote the Gauss-sum over $\bbbf_q$
with multiplicative character $\omega^m$. Then, for any integer $m$ we have
$$g_q(m)=-\prod_{i=0}^{f-1}\pi^{(p-1)\left\{{p^im\over q-1}\right\}}
\Gamma_p\left(\left\{{p^im\over q-1}\right\}\right).$$
Here $\{x\}=x-\lfloor x\rfloor$ is the fractional part of $x$.
In particular, when $q=p$ we get
$$g_p(m)=-\pi^{(p-1)\left\{{m\over p-1}\right\}}
\Gamma_p\left(\left\{{m\over p-1}\right\}\right).$$
\end{theorem1}

See Henri Cohen's book \cite{cohen2} for a proof.
When $p$ does not divide the common denominator of the $\alpha_i,\beta_j$ one easily
writes down a $p$-adic version of our hypergeometric sum for the case $q=p$.

\begin{definition1}\label{padic}
We define $G_p(\boldalpha,\boldbeta|t)$ by the sum
$${1\over 1-p}\sum_{m=0}^{q-2}\omega((-1)^dt)^m
(-p)^{\Lambda(m)}\prod_{i=1}^d{\Gamma_p\left(\left\{\alpha_i+{m\over p-1}\right\}\right)
\over\Gamma_p(\{\alpha_i\})}
{\Gamma_p\left(\left\{-\beta_i-{m\over p-1}\right\}\right)\over
\Gamma_p(\{-\beta_i\})},$$
where 
$$\Lambda(m)=\sum_{i=1}^d
\left\{\alpha_i+{m\over p-1}\right\}-\{\alpha_i\}
+\left\{-\beta_i-{m\over p-1}\right\}-\{-\beta_i\}.$$
\end{definition1}

Note that 
$$\Lambda(m)=\sum_{i=1}^d
-\left\lfloor\alpha_i+{m\over p-1}\right\rfloor+\floor{\alpha_i}
-\left\lfloor-\beta_i-{m\over p-1}\right\rfloor+\floor{-\beta_i}.$$
In particular $\Lambda(m)\in\bbbz$.
Definition \ref{padic} almost coincides with McCarthy's function
$\,_dG_d$ from \cite[Def 1.1]{mccarthy}
in the sense that our function coincides with $\,_dG_d(1/t)$. We prefer to adhere to the 
definition given above. The advantage of Definition \ref{padic} is that Assumption
\ref{divisionassumption} is not required, it is well-defined for
all parameters $\alpha_i,\beta_j$ as long as they are $p$-adic integers. Define
$$\delta=\delta(\boldalpha,\boldbeta)=
\max_{x\in[0,1]}\sum_{i=1}^d\floor{x+\alpha_i}-\floor{\alpha_i}
+\floor{-x-\beta_i}-\floor{-\beta_i}.$$
Then, using Definition \ref{padic} and the fact that $-\Lambda(m)\le\delta$ one easily
deduces that $p^{\delta}G_p(\boldalpha,\boldbeta|t)$ is a $p$-adic integer.
In \cite[Prop 3.1]{mccarthy} we find this in a slightly different formulation.
However, it is not clear from the definition whether this value is algebraic or not over
$\bbbq$. It is the purpose of the present note to be a bit more
specific by proving the following theorem.

\begin{theorem1}\label{main}
Let notations be as above and let $K$ be the field extension of $\bbbq$ generated by
the coefficients of $A(x)$ and $B(x)$. Suppose $p$ splits in $K$, i.e. $p$ factors into
$[K:\bbbq]$ distinct prime ideals in the integers of $K$.
Let $\Delta=\max_k\delta(k\boldalpha,k\boldbeta)$ over all integers $k$ relatively
prime with the common denominator of the $\alpha_i,\beta_j$.
Then $p^{\Delta}G_p(\boldalpha,\boldbeta|t)$ is an algebraic integer in $K$. 
\end{theorem1}

For the proof we construct in Section \ref{algebras} a generalization of the
hypergeometric function $H_q(A,B|t)$ involving two semisimple finite algebras
$A$ and $B$ over $\bbbf_q$. We show that it belongs to $K$ and then, in Section
\ref{mainproof} identify its $p$-adic evaluation with $G_p(\boldalpha,\boldbeta|t)$.

\section{Gauss sums on finite algebras}\label{algebras}
The main idea of the proof of Theorem \ref{main}
is to use Gauss sums on finite commutative algebras over $\bbbf_q$ with $1$.
Let $A$ be such an algebra. For any $x\in A$ we define
the trace $\tr(x)$ and norm $N(x)$ as the trace and
norm of the $\bbbf_p$-linear map given by multiplication with $x$ on $A$.

Choose an additive character $\psi$ on $A$ which is {\it primitive}. That is,
to any ideal $I\subset A, I\ne(0)$ there exists $x\in I$ such that $\psi(x)\ne1$.
Any other non-degenerate additive character is of the form $\psi(ax)$ with $a\in A^{\times}$. 
A multiplicative character $\chi$ is called {\it primitive} if its
kernel does not contain any subgroup of the form $\{1+a|a\in I\}$ for some
non-zero ideal $I$ in $A$. 

For any multiplicative character $\chi$ on $A^{\times}$ we can define a 
Gauss sum
$$g_A(\psi,\chi)=\sum_{x\in A^{\times}}\psi(x)\chi(x).$$
When $A$ is not semisimple, the Gauss sum can be $0$, as illustrated by the
following example.

\begin{example1}
Let $A=\bbbf_p[x]/(x^2)$. Choose the additive character $\psi(a+bx)=\zeta_p^b$.
It is easy to see that this is
a primitive character. Note that $a+bx\in A^{\times}\iff
a\in\bbbf_p^{\times}$. Let $\chi$ be a nontrivial multiplicative
character on $\bbbf_p^{\times}$ and extend it to $A^{\times}$ by $\chi(a+bx)=\chi(a)$. 
Then
$$g_A(\psi,\chi)=\sum_{a\in\bbbf_p^{\times},b\in\bbbf_p}\zeta_p^b\chi(a)=0.$$
\end{example1}

So we restrict ourselves to semisimple algebras. These are precisely the finite
sums of finite field extensions of $\bbbf_q$. In this case there is an obvious
choice for the additive character.

\begin{lemma1}
Suppose $A$ is a direct sum of finite field extensions of $\bbbf_q$.
Then $\psi(x)=\zeta_p^{\tr(x)}$ is a primitive additive character.
\end{lemma1}

\begin{proof1}{}
Let $A\cong \oplus_{i=1}^r F_i$ with $F_i$ a finite field extension of $\bbbf_q$
for all $i$. Then $\psi(x)=\zeta_p^{\tr_1(x_1)+\cdots+\tr_r(x_r)}$, where
$\tr_i$ stands for the trace function on $F_i$. If $\psi$ were not
primitive then there exists $a\in A,a\ne0$ such that $\psi(ax)=1$ for all $x\in A$.
Suppose $a=(a_1,\ldots,a_r)$ and assume, without loss of generality, $a_1\ne0$.
Then $\psi(x,0,\ldots,0)=\zeta_p^{\tr(a_1x)}=1$ for all $x\in F_1$. By the properties
of the trace of a field this is not possible. 
\end{proof1}

From now on we use the trace character on a semisimple algebra $A$ as additive character
and write $g_A(\chi)$ for the Gauss sum. So we dropped the dependence of the Gauss sum
on the aditive character. The only amount of freedom in the additive character rests on
the choice of $\zeta_p$.

\begin{proposition1}\label{normgausssum}
Let $A$ be a direct sum of finite fields over $\bbbf_q$ and $\psi(x)=\zeta_p^{\tr(x)}$
the additive character. Let $\chi$ a multiplicative character. Then there
exists a non-negative integer $f$ such that
$$|g_A(\chi)|^2=q^f.$$
\end{proposition1}

\begin{proof1}{}
Again, write $A=\oplus_{i=1}^r F_i$. Then $\chi$ can be written as
$\chi(x_1,\ldots,x_r)=\chi_1(x_1)\cdots\chi_r(x_r)$, where $\chi_i$
is a multplicative character on $F_i^{\times}$. This implies that
$$g_A(\psi,\chi)=\prod_{i=1}^r g(\chi_i),$$
where $g(\chi_i)$ is the usual Gauss sum on the field $F_i$.
The additive character on $F_i$ is $\zeta_p^{\tr_i(x)}$ with the same
choice of $\zeta_p$ for each $i$. 
Our assertion follows directly. 
\end{proof1}

Choose two finite semisimple algebras $A,B$ over $\bbbf_q$. Choose the trace characters
on each of them with the same choice of $\zeta_p$ and call them $\psi_A,\psi_B$.
Let $\chi_A,\chi_B$ be multiplicative characters on $A^{\times},B^{\times}$. 
Denote the norms on $A,B$ by $N_A,N_B$.
 
\begin{definition1}\label{algebraHGF}
We define
$$H_q(A,B|t)={-1\over g_A(\chi_A)g_B(\chi_B)}\sum_{x\in A^{\times},y\in B^{\times},tN_A(x)=N_B(y)}
\psi_A(x)\psi_B(-y)\chi_A(x)\overline{\chi_B(y)},$$
for any $t\in\bbbf_q$.
\end{definition1}

The following theorem gives its Fourier expansion in $t$.

\begin{theorem1}
Let $\omega$ be a generator of the multiplicative characters on $\fqstar$.
When the context is clear we denote both functions $\omega(N_A(x))$ and $\omega(N_B(y))$
by $\omega_N$. We then have,
$$H_q(A,B|t)={1\over 1-q}\sum_{m=0}^{q-2}{g_A(\chi_A\omega_N^m)g_B(\overline{\chi_B}\omega_N^{-m})
\over g_A(\chi_A)g_B(\chi_B)}\omega(N_B(-1)t)^m.$$
\end{theorem1}

\begin{proof1}{}
We compute the Fourier expansion $\sum_{m=0}^{q-2}c_m\omega(t)^m$ of $H_q(A,B|t)$.
The coefficient $c_m$ can be computed using 
$$c_m={1\over q-1}\sum_{t\in\fqstar}H_q(A,B|t)\omega(t)^{-m}.$$
When we substitute the definition for $H_q(A,B|t)$ in the summation over $t$,
we get a summation over $t\in\fqstar,x\in A^{\times},y\in B^{\times}$ with the
restriction $tN_A(x)=N_B(y)$. So we might as well substitute $t=N_B(y)/N_A(x)$ and sum
over $x,y$. We get,
$$
c_m={1\over 1-q}\sum_{x\in A^{\times},y\in B^{\times}}{1\over g_Ag_B}
\psi_A(x)\psi_B(-y)\chi_A(x)\chi_B(y)^{-1}\omega(N_A(x))^m\omega(N_B(y))^{-m}.$$
The summation over $x$ yields $g_A(\chi_A\omega_N^m)$. To sum over $y$ we first
replace $y$ by $-y$ and then perform the summation. We get 
$\omega(N_B(-1))^mg_B(\overline{\chi_B}\omega_N^{-m})$. This proves our theorem.
\end{proof1}

\begin{example1}\label{example}
As in the previous section take two multisets of hypergeometric parameters $\boldalpha,\boldbeta$.
Suppose that $(q-1)\alpha_i,(q-1)\beta_j$ are in $\bbbz$ for all $i,j$.
Take $A=B=\bbbf_q^d$, the direct sum of $d$ copies of $\bbbf_q$ with componentwise
addition and multiplication. The norm on $A,B$ is given by $N(x_1,\ldots,x_d)=
x_1\cdots x_d$. In particular $N_B(-1)=(-1)^d$.
For both $A,B$ we take the additive character
$\psi(x_1,\ldots,x_d)=\zeta_p^{\tr(x_1+\cdots+x_d)}$, where $\tr$ the trace
function on $\bbbf_q$.
As multiplicative characters we take 
$$\chi_A(x_1,\ldots,x_d)=\prod_{i=1}^d\omega(x_i)^{(q-1)\alpha_i},\qquad
\chi_B(x_1,\ldots,x_d)=\prod_{j=1}^d\omega(y_j)^{(q-1)\beta_j}.$$

An easy calculation shows that $g_A(\chi_A\omega_N^m)=\prod_{i=1}^dg(m+(q-1)\alpha_i)$
and similarly for $g_B$. So we see that we recover the finite hypergeometric sum of
the previous section.
\end{example1}

\begin{lemma1}
Suppose $\dim_{\bbbf_q}(A)=\dim_{\bbbf_q}(B)$. Then $H_q(A,B|t)$ does not depend on the 
choice of $\zeta_p$ in the additive characters.

As a corollary, in this equi-dimensional case the values of $H_q(A,B|t)$ 
are contained in the field generated by the charactervalues of $\chi_A,\chi_B$.
\end{lemma1}

\begin{proof1}{}
When we choose $\zeta_p^a, a\in\bbbf_p^{\times}$ instead of $\zeta_p$
in the definition of the additive character
it is easy to check that $g_A(\chi_A)$ gets replaced by $\chi_A(a)^{-1}g_A(\chi_A)$.
And similarly for $B$.
As a corollary any term in the sum in the hypergeometric sum in Theorem \ref{algebraHGF}
is multiplied by $\omega(N_B(a)/N_A(a))^m$.
Since $a\in\bbbf_p$ is a scalar, $N_A(a)=N_B(a)=a^d$, where
$d=\dim_{\bbbf_q}(A)=\dim_{\bbbf_q}(B)$. Hence, in the case of equal dimensions
of $A,B$ the multiplication factor is $1$. 

Let $\sigma\in{\rm Gal}(\overline{\bbbq}/\bbbq)$ be such that it fixes the values of
$\chi_A,\chi_B$ but sends $\zeta_p$ to $\zeta_p^{a}$.
According to the above calculation $H_q(A,B|t)$ is fixed
under this substitution and hence under $\sigma$. 
\end{proof1} 

Let us return momentarily to Example \ref{example} and suppose that the parameters
$\boldalpha$ have the property that $k\boldalpha\is\boldalpha\mod{\bbbz},
k\boldbeta\is\boldbeta\mod{\bbbz}$ for all $k$ relative prime with the 
common denominator of the $\alpha_i,\beta_j$.
Then, for any $\sigma\in\galQ$ there exists a
permutation $\rho$ of the summands of $A=\oplus_{i=1}^d\bbbf_p$ such that
$\chi_A(\rho(x))=\chi_A(x)^{\sigma}$ for all $x\in A^{\times}$. A similar 
permutation exists for $B$. Notice also that $\tr(\rho(x))=\tr(x)$ 
and $N(\rho(x))=N(x)$.  

A similar situation arises in the case $A=\bbbf_{p^r}$ as $\bbbf_p$-algebra.
Let $\chi_A$ be a character of order $d$ dividing $p^r-1$. Let $\rho$ be
the $p$-th power Frobenius on $A$, then $\chi_A(\rho(x))=\chi_A(x)^p$,
a conjugate of $\chi_A(x)$ for all $x\in A^{\times}$. Notice also that
$\tr(\rho(x))=x$ and $N(\rho(x))=N(x)$.

\begin{definition1}
Let $A$ be a finite dimensional $\bbbf_q$-algebra. A ring automorphism $\rho:A\to A$
is called an $\bbbf_q$-automorphism if it is $\bbbf_q$-linear and it fixes both
norm and trace of $A$. 
\end{definition1}

\begin{proposition1}\label{fieldofdefinition}
Let $A,B$ be finite commutative semisimple $\bbbf_q$-algebras. Let $\chi_A,\chi_B$ be 
multiplicative characters. Consider the subgroup $G$ of $\galQ$ of elements
$\sigma$ for which there exists an $\bbbf_q$-automorphisms $\rho_A$ of $A$ 
and $\rho_B$ of $B$ with the property that $\chi_A(\rho_A(x))=\chi_A(x)^{\sigma}$
and $\chi_B(\rho_B(x))=\chi_B(x)^{\sigma}$ for every $\sigma\in G$. Then
$H_q(A,B|t)$ lies in the fixed field of $G$ for every $t\in\bbbf_q^{\times}$. 
\end{proposition1}

\begin{proof1}{}
Let $\sigma\in G$. We first compute the action of $\sigma$ on $g_A(\chi_A)$. 
Suppose that $\sigma(\zeta_p)=\zeta_p^a$.
\begin{eqnarray*}
g_A(\chi_A)^{\sigma}&=&\sum_{x\in A^{\times}}\zeta_p^{a\tr(x)}\chi_A(x)^{\sigma}\\
&=&\sum_{x\in A^{\times}}\zeta_p^{a\tr(x)}\chi_A(\rho(x))\\
&=&\sum_{x\in A^{\times}}\zeta_p^{\tr(\rho^{-1}(x))}\chi_A(a^{-1}x)\\
&=&\chi_A(a)^{-1}g_A(\chi_A)
\end{eqnarray*}
A similar calculation holds for $B$. Now apply $\sigma$ to the terms in the
sum in Definition \ref{algebraHGF}. A similar calculation as above shows
that the sum gets multiplied with $\chi_A(a)^{-1}\chi_B(a)^{-1}$. This 
cancels the factor coming from $g_A(\chi_A)g_B(\chi_B)$. Hence $H_q(A,B|t)$
is fixed under all $\sigma\in G$.
\end{proof1}

\section{Proof of Theorem \ref{main}}\label{mainproof}
We use the notations from the introduction. In particular
$$A(x)=\prod_{j=1}^d(x-e^{2\pi i\alpha_j}),\quad B(x)=\prod_{j=1}^d
(x-e^{2\pi i\beta_j})$$
and $K$ is the field generated by the coefficients of $A(x)$ and $B(x)$.
Let $p$ be a prime which splits completely in $K$. Then
we can consider $A(x)$ as element of $\bbbf_p[x]$. 
Let $A(x)=A_1(x)\cdots A_r(x)$ be the irreducible factorization of $A(x)$
in $\bbbf_p[x]$. For the $\bbbf_p$-algebra we take
$\oplus_{i=1}^r\bbbf_p[x]/(A_i(x))$. The construction of a multiplicative
character on $A$ is as follows. First we choose a multiplicative character
$\omega$ on $\overline{\bbbf_p}$ such that its restriction to $\bbbf_{p^r}$
has order $p^r-1$ for all $r\ge1$ and fix in the remainder of the proof.

Since $p$ splits in $K$ multiplication by $p$ gives a permutation of the 
multiset $\boldalpha$ modulo $\bbbz$. Under this action $\boldalpha\mod{\bbbz}$ decomposes
into a union of orbits, which we call $p$-orbits. Let $O$ be such a $p$-orbit. Then 
$\prod_{\alpha\in O}(x-e^{2\pi i\alpha})$ is a polynomial and $p$ splits in the field
generated by its coefficients. So we can consider it modulo a prime ideal dividing $p$ and
hence as an element of $\bbbf_p[x]$. It is one of the factors $A_i(x)$ of the mod $p$
factorization of $A(x)$. The orbit $O$ will now be denoted by $O_i$. There are $r$
orbits and we renumber the indices of the $\alpha_i$ such that $\alpha_i\in O_i$
for $i=1,\ldots,r$.
On $\bbbf_p[x]/(A_i)$ we define the multiplicative
character $\chi_i=\omega^{\alpha_i(q_i-1)}$, where $q_i=p^{\deg(A_i)}$.
If we would have chosen $p\alpha_i$ instead of $\alpha_i$, the new character
would simply consist of the Frobenius transform followed by $\chi_i$. 
For the character $\chi_A$ on $A=\sum_{i=1}^r\bbbf_p[x]/(A_i)$ we choose 
$$\chi_A(x_1,\ldots,x_r)=\prod_{i=1}^r\omega(x_i)^{\alpha_i(q_i-1)}.$$

Let $\sigma\in\galK$. It acts as $\omega(x)\mapsto
\omega(x)^k$ for some integer $k$. Hence
$$\chi_A(x_1,\ldots,x_r)^{\sigma}=\prod_{i=1}^r\omega(x_i)^{k\alpha_i(q_i-1)}.$$
This permutes the factors by a permutation $s\in S_r$ and we get
$$\chi_A(x_1,\ldots,x_r)^{\sigma}=\prod_{i=1}^r\omega\left(x_{s^{-1}(i)}\right)^{p^{l_i}\alpha_i
(q_i-1)},$$
where $0\le l_i<\deg(A_i)$ for each $i$. We used $q_{s(i)}=q_i$. We finally get
$$\chi_A(x_1,\ldots,x_r)^{\sigma}=
\prod_{i=1}^r\omega\left(x_{s^{-1}(i)}^{p^{l_i}}\right)^{\alpha_i(q_i-1)}
=\chi_A\left(x_{s^{-1}(1)}^{p^{l_i}},\ldots,x_{s^{-1}(r)}^{p^{l_r}}\right).$$
In other words, $\chi_A(x)^{\sigma}=\chi_A(\rho(x))$ for a suitable $\bbbf_p$-automorphism
$\rho$ of $A$. Notice that norm and trace of $A$ are preserved by $\rho$.
A similar construction can be performed for $B(x)$. 
According to Proposition \ref{fieldofdefinition} we get $H_p(A,B|t)\in K$
for all $t\in\fqstar$. 

In order to connect to the $p$-adic function $G_p$ we take the inverse of the
Teichm\"uller character for $\omega$ and compute the terms given in 
Definition \ref{algebraHGF} $p$-adically. The Gauss sum $g_A(\chi_A\omega_N^m)$
is the product of ordinary Gauss sums of the form $g(\omega^{(q-1)\alpha+m(1+p+\cdots+p^{l-1})})$
over the field $\bbbf_q$ with $q=p^l$. The occurrence of $m(1+p+\cdots+p^{l-1})$ is due
to $\omega(N_{\bbbf_q/\bbbf_p}(x)^m)=\omega(x)^{m(1+\cdots+p^{l-1})}$. The Gross-Koblitz
theorem for Gauss sums over $\bbbf_q$ with $q=p^l$ gives us
$$g_q(\omega^a)=-\prod_{i=0}^{l-1}\pi^{\left\{{p^ia\over q-1}\right\}}
\Gamma_p\left(\left\{{p^ia\over q-1}\right\}\right)$$
for every integer $a$. 
When applied to $a=(q-1)\alpha+m(q-1)/(p-1)$ this amounts to
$$-\prod_{i=0}^{l-1}\pi^{\left\{p^i\alpha+{m\over p-1}\right\}}
\Gamma_p\left(\left\{p^i\alpha+{m\over p-1}\right\}\right).$$
Note that this is a product over the $p$-orbit containing $\alpha$ and each
factor is precisely of the type that occur in the definition of the $p$-adic
hypergeometric sum. A similar story goes for $B(x)$. As a result we get
$${g_A(\chi_A\omega_N^m)g_B(\overline{\chi_B}\omega_N^{-m})
\over g_A(\chi_A)g_B(\overline{\chi_B})}=(-p)^{\Lambda(m)}
\prod_{i=1}^d{\Gamma_p\left(\left\{\alpha_i+{m\over p-1}\right\}\right)
\Gamma_p\left(\left\{-\beta_i-{m\over p-1}\right\}\right)\over
\Gamma_p\left(\left\{\alpha_i\right\}\right)\Gamma_p\left(\left\{-\beta_i\right\}\right)},$$
where $\Lambda(m)$ is as defined in the introduction. 
So we find that $p$-adically
$$H_p(A,B|t)=G_p(\boldalpha,\boldbeta|t).$$
Hence we conclude that the values of $G_p$ are in $K$. It remains to give an estimate
for the denominator. From Definition \ref{algebraHGF} it follows that $H_p(A,B|t)$
has the denominator
$g_Ag_B$. Hence, by Proposition \ref{normgausssum} there exists a power
$p^r$ of $p$ such $p^rH_p(A,B|t)$ is an integer
in $K$. It remains to determine a value for $r$.
The conjugates of $H_p(A,B|t)$ are obtained by taking $\chi_A^k,\chi_B^k$
as multiplicative characters. The corresponding hypergeometric parameters are
$k\boldalpha,k\boldbeta$. From McCarthy's work it follows that 
$p^{\Delta}G_p(k\boldalpha,k\boldbeta|t)$ is a $p$ adic integer for all
$k$ relatively prime to the common denominator of $\alpha_i,\beta_j$.
This implies that $p^{\Delta}H_p(A,B|t)$ is an algebraic integer in $K$.

\end{document}